\documentclass[a4paper]{amsart}

\title[BRST cohomologies for symplectic reflection algebras]{%
BRST cohomologies for symplectic reflection algebras and quantizations of hypertoric varieties}
\author{Toshiro Kuwabara}
\thanks{The author was partially supported by 
Grant-in-Aid for Young Scientist (B) 21740013,
Japan Society for the Promotion of Science}
\thanks{The author was partially supported by 
GCOE Program `Fostering top leaders in mathematics', Kyoto University.
}
\thanks{This work was
 partially supported by Basic Science Research Program through the National Research Foundation of Korea (NRF) grant funded by the Korea government (MEST)(2011-0027952).}
\address{Department of Mathematics, National Research University -- Higher School
of Economics, 20 Myasnitskaya st., Moscow 101000, Russia.}
\email{toshiro.kuwa@gmail.com}

\usepackage{latexsym}
\usepackage{amsmath}
\usepackage{amsfonts}
\usepackage{amssymb}
\usepackage{amsxtra}
\usepackage{mathrsfs}
\usepackage [matrix,arrow]{xy}
\usepackage{pdfsync}

\newtheorem{definition}{Definition}[section]
\newtheorem{proposition}[definition]{Proposition}
\newtheorem{theorem}[definition]{Theorem}

\newtheorem{lemma}[definition]{Lemma}
\newtheorem{remark}[definition]{Remark}

\newcommand{\refprop}[1]{Proposition~\ref{#1}}
\newcommand{\refthm}[1]{Theorem~\ref{#1}}

\newcommand{\reflemma}[1]{Lemma~\ref{#1}}

\newcommand{\refeq}[1]{(\ref{#1})}

\newcommand{\refsec}[1]{Section~\ref{#1}}

\newcommand{\W}{\scW}
\newcommand{\pZ}{{}'Z}
\newcommand{\pB}{{}'B}
\newcommand{\pE}{{}'E}

\newcommand{\RHom}[1]{\mathrm{R^{#1}Hom}}
\newcommand{\LTensor}{\stackrel{L}{\otimes}}
\newcommand{\blkbar}{\raisebox{0.5ex}{\rule{2ex}{0.4pt}}}

\newcommand{\defeq}{\ensuremath{\underset{\mathrm{def}}{=}}}
\newcommand{\C}{{\mathbb C}}

\newcommand{\Q}{{\mathbb Q}}

\newcommand{\Z}{{\mathbb Z}}

\newcommand{\bbG}{{\mathbb G}}

\newcommand{\bbX}{{\mathbb X}}

\newcommand{\bfv}{\mathbf{v}}
\newcommand{\bfw}{\mathbf{w}}

\newcommand{\calC}{\mathcal{C}}
\newcommand{\calD}{\mathcal{D}}

\newcommand{\calF}{\mathcal{F}}

\newcommand{\calH}{\mathcal{H}}

\newcommand{\calO}{\mathcal{O}}

\newcommand{\calW}{\mathcal{W}}

\newcommand{\frD}{\mathfrak{D}}

\newcommand{\frL}{\mathfrak{L}}

\newcommand{\frS}{\mathfrak{S}}

\newcommand{\frU}{\mathfrak{U}}

\newcommand{\frX}{\mathfrak{X}}

\newcommand{\frg}{\mathfrak{g}}

\newcommand{\gl}{\mathfrak{gl}}

\newcommand{\scF}{\mathscr{F}}

\newcommand{\scL}{\mathscr{L}}
\newcommand{\scM}{\mathscr{M}}
\newcommand{\scN}{\mathscr{N}}

\newcommand{\scW}{\mathscr{W}}

\DeclareMathOperator{\Tr}{Tr}
\DeclareMathOperator{\gEnd}{End}
\DeclareMathOperator{\lEnd}{\operatorname{\mathscr{E}\kern-.1pc\mathit{nd}}}
\DeclareMathOperator{\gHom}{Hom}
\DeclareMathOperator{\lHom}{\operatorname{\mathscr{H}\kern-.1pc\mathit{om}}}
\DeclareMathOperator{\Mod}{Mod}
\newcommand{\mmod}{\text{-}\mathrm{mod}}
\DeclareMathOperator{\Ad}{Ad}
\DeclareMathOperator{\ad}{ad}
\DeclareMathOperator{\Aut}{Aut}
\DeclareMathOperator{\Der}{Der}

\DeclareMathOperator{\Spec}{Spec}
\DeclareMathOperator{\Proj}{Proj}
\DeclareMathOperator{\Lie}{Lie}

\DeclareMathOperator{\Dim}{dim}
\DeclareMathOperator{\Ker}{Ker}
\renewcommand{\Im}{\operatorname{Im}}
\DeclareMathOperator{\Span}{Span}
\DeclareMathOperator{\gr}{gr}

\newcommand{\opp}{\mathrm{opp}}

\newcommand{\isoto}[1][]{\mathop{\xrightarrow[#1]%
{\rule{0pt}{.9ex}%
{\raisebox{-.4ex}[0ex][-.6ex]{$\mspace{3mu}\sim\mspace{3mu}$}}}}}
\newcommand{\scbul}{{\,\raise1pt\hbox{$\scriptscriptstyle\bullet$}\,}}
\newcommand{\prolim}[1][]{\mathop{\varprojlim}\limits_{#1}}

\begin{document}
\begin{abstract}
 We study algebras constructed by quantum Hamiltonian reduction
 associated with symplectic quotients of symplectic vector spaces,
 including deformed preprojective algebras, symplectic reflection
 algebras (rational Cherednik algebras),
 and quantization of hypertoric varieties introduced by Musson and Van den Bergh
 in \cite{MVdB}.
 We determine BRST cohomologies associated with
 these quantum Hamiltonian reductions. To compute these BRST cohomologies,
 we make use of method of deformation quantization (DQ-algebras) and F-action
 studied by Kashiwara and Rouquier in \cite{KR}, and Gordon and Losev in \cite{GL}.
\end{abstract}

\maketitle

\section{Introduction}
\label{sec:introduction}

Symplectic reflection algebras were introduced by Etingof and Ginzburg
in \cite{EG}. For a symplectic vector space $V$ and a finite group
$W$ generated by symplectic reflections of $V$, they defined a symplectic
reflection algebra as noncommutative deformation of symplectic quotient
singularity $V / W$. For $W = G(\ell, 1, n)$, a complex reflection group,
the associated symplectic reflection algebra is also called a rational
Cherednik algebra.

In \cite{GG2}, \cite{Go2} and \cite{EGGO}, 
they studied construction of the symplectic reflection
algebra associated to $W=\Gamma \wr \frS_n$, the wreath product of
a finite subgroup $\Gamma \subset SL_2(\C)$ with the symmetric group $\frS_n$,
 by quantum Hamiltonian reduction.
These results were generalization
of a result of \cite{H} which studied construction of deformed 
preprojective algebras by quantum Hamiltonian reduction.

Hypertoric varieties are quaternionic analogue of toric varieties and
are constructed by Hamiltonian reduction of a symplectic vector space
with an action of a torus.
Musson and Van den Bergh constructed quantization of the hypertoric 
varieties by using quantum Hamiltonian reduction in \cite{MVdB}.

We briefly review these quantum Hamiltonian reductions.
Let $V$ be a vector space over $\C$ with an action of a reductive algebraic
group $G$. Let $\frD(V)$ be the algebra of algebraic differential operators
on $V$. Then the action of $G$ on $V$ induces an action of $G$ on $\frD(V)$.
Moreover, we have a Lie algebra homomorphism 
$\mu_{\frD} : \frg (\defeq Lie G) \longrightarrow \frD(V)$ which we call
a quantized moment map. We consider the subalgebra of $G$-invariant elements
in the quotient of $\frD(V)$ by the image of $\mu_{\frD}$,
\[
 \frD(X_0, c) = \bigl(\frD(V) / \frD(V) (\mu_{\frD} + c)(\frg)\bigr)^G
\]
where $c : \frg \longrightarrow \C$ is a character of $\frg$. This algebra
can be regarded as quantization of a Hamiltonian reduction 
$X_0 = \mu_{T^*V}^{-1}(0) /\!/ G$ of the cotangent bundle $T^*V$, where  
$\mu_{T^*V}: T^* V \longrightarrow \frg^*$ is the corresponding moment map. 
Thus, we call the algebra $\frD(X_0, c)$ a quantum Hamiltonian reduction
of $\frD(V)$ with respect to the $G$-action.

We study BRST cohomologies associated to the quantum Hamiltonian reduction.
The BRST cohomologies were first introduced by theoretical physicists and now play an important
role in mathematical physics and representation theory to give a
cohomological description of quantum Hamiltonian reduction.
(see e.g. \cite{KS}).

For finite $\calW$-algebras, quantizations of Slodowy slices, vanishing
of Lie algebra homologies and cohomologies associated with
quantum Hamiltonian reduction corresponding to them was proved by Gan and Ginzburg in \cite{GG}
and \cite{G}. These results immediately imply that the BRST cohomologies associated
with the finite $\calW$-algebras are concentrated in degree $0$. 

Motivated by their results, we study BRST cohomologies associated with
the quantum Hamiltonian reductions corresponding to the deformed preprojective
algebras, the symplectic reflection algebras and quantized hypertoric algebras.
For these algebras, positive BRST cohomologies do not vanish in contrast
to ones for the finite W-algebras. Nevertheless, we prove vanishing of 
negative BRST cohomologies for these algebras and determine positive
BRST cohomologies explicitly. 

Recently, (micro-)localization of these quantum Hamiltonian reductions was
studied by using deformation quantization (DQ-algebras)
based on the ideas of \cite{KR} and \cite{Losev} (see also \cite{DK}, \cite{BK} and \cite{BLPW}).
A DQ-algebra is a sheaf of noncommutative $\C((\hbar))$-algebras on a (symplectic)
resolution $X$ of the singularities of $X_0$ which is isomorphic to
$\calO_{X} \otimes_{\C} \C((\hbar))$ as a sheaf of vector spaces. 
By considering analogue of the quantum Hamiltonian reduction for $\frD(X_0, c)$, we
can construct a DQ-algebra $\scW_{X,c}$ on $X$.
Moreover, by introducing an equivariant $\C^*$-action on the DQ-algebra,
we can obtain the algebra $\frD(X_0, c)$ as the algebra of its $\C^*$-invariant
global sections. 

We briefly review the construction. On the symplectic
vector space $T^*V$, we have a canonical DQ-algebra $\scW_{T^*V}$ with
an equivariant $G$-action.
Consider 
an open subset $\frX$ consisting of semistable points with respect to the
$G$-action. The symplectic manifold $X$ can be constructed as
$X = (\mu_{T^*V}^{-1}(0) \cap \frX) / G$. 
Let $p : \mu_{T^*V}^{-1}(0) \cap \frX \longrightarrow X$ be the projection.
Set $\scW_{\frX} \defeq \scW_{T^*V} \vert_{\frX}$. We have a homomorphism of
algebras $\frD(V) \longrightarrow \scW_{T^*V}(T^*V)$ and this induces a
quantized moment map $\mu_{\scW} : \frg \longrightarrow \scW_{\frX}(\frX)$.
For the character $c$ of $\frg$,
the DQ-algebra $\scW_{X,c}$ on $X$ is a sheaf
associated to the presheaf 
\[
 U \mapsto \scW_{X, c}(U) = \bigl(\scW_{\frX}(\frU) / \scW_{\frX}(\frU) (\mu_{\scW} + c)(\frg)\bigr)^G
\]
where $\frU \subset \frX$ is an open subset such that $p^{-1}(U) = \frU \cap \mu^{-1}_{T^*V}(0)$. Under some geometrical conditions, we obtain an isomorphism
$\scW_{X,c}(X)^{\C^*} \simeq \frD(X_0, c)$.

Associating with the above two quantum Hamiltonian reductions, we define
BRST cohomologies $H_{BRST, c}^{\bullet}(\frg, \frD(V))$ and
$\calH_{BRST, c}^{\bullet}(\frg, \scW_{\frX})$ (see \refsec{sec:brst-cohomology}
for the definition of the BRST cohomologies).
The main result of the paper is the following two theorems.

\begin{theorem}[\refthm{thm:local-BRST}]
 We have the following isomorphism of sheaves on $X$,
 \[
  \calH^n_{BRST,c}(\frg, \scW_{\frX})
 \simeq \scW_{X,c} \otimes_{\C} H_{DR}^n(G)
 \]
 where $H_{DR}^{\bullet}(G)$ is the (algebraic) de Rham cohomology of $G$.
\end{theorem}

By using the two kinds of equivalences between the category of 
$\C^*$-equivariant $\scW_{X,c}$-modules and the category of $\frD(X_0, c)$-modules, which
are called abelian and derived $\scW$-affinities, we also have the following
explicit description of the BRST cohomology $H_{BRST,c}^{\bullet}(\frg, \frD(V))$.

\begin{theorem}[\refthm{thm:global-BRST}]
 We have the following isomorphism of $\C$-algebras,
 \[
  H^{\bullet}_{BRST,c}(\frg, \frD(V))
 \simeq \frD(X_0, c) \otimes_{\C} H_{DR}^{\bullet}(G).
 \]
\end{theorem}

This paper is organized as follows: In \refsec{sec:w-algebras}, we review the definition
and basic properties of DQ-algebras. In \refsec{sec:setting-q-Hamilton-red}, we introduce
the quantum Hamiltonian reduction both for usual algebras of differential operators and
for DQ-algebras. \refsec{sec:alg-q-Hamilton} is review of construction of 
symplectic reflection algebras and quantized hypertoric algebras by the quantum Hamiltonian
reduction. We do {\em not} use facts in this section later for proving our main theorem,
and thus, the reader can skip this section. In \refsec{sec:brst-cohomology}, we introduce
BRST cohomologies associated with the quantum Hamiltonian reduction. 
\refsec{sec:comp-brst-cohom} is the main part of this paper. In \refsec{sec:vanish-neg-cohom},
we prove that negative BRST cohomologies vanish if the moment map is flat, and in
\refsec{sec:pos-BRST-cohom}, we determine positive BRST cohomologies.
Finally, we apply these results for the algebras, which reviewed in \refsec{sec:alg-q-Hamilton},
and determine the BRST cohomologies associated to them explicitly in \refsec{sec:examples}.

\subsection*{Acknowledgments}

The author is deeply grateful to Tomoyuki Arakawa for 
his introduction to the BRST cohomologies and valuable 
discussions and comments. He also thanks Gwyn Bellamy for his valuable comments and his advice about
geometric facts for the quiver varieties and hypertoric varieties.
He is also grateful to Masaki Kashiwara and Ivan Losev for valuable discussions. 
The author is grateful to Max Planck Institute for Mathematics for hospitality during the
period from February to May 2011. 

\section{Preliminaries}

Let $G$ be a group and $V$ be a $G$-module. We denote the subset of
all $G$-invariant elements of $V$ by $V^G$. For a character 
$\theta : G \longrightarrow \C^*$, we denote the subset of all
$G$-semi-invariant element belonging to the character $\theta$ by
$V^{G, \theta}$. For an element $v \in V$, let 
$G_v = \{ g \in G \;\vert\; g \cdot v = v \}$ be the stabilizer of $v$.

For a Lie algebra $\frg$, we denote the universal enveloping algebra of $\frg$
by $U(\frg)$. We also denote the symmetric algebra over $\frg$ 
by $S(\frg) = \C[\frg^*]$ and the exterior algebra over $\frg$ by $\Lambda(\frg)$.

For a commutative algebra $A$ over $\C$, let $\Spec A$ be the affine 
scheme associated to $A$. For a commutative graded algebra 
$A = \bigoplus_{n \in \Z_{\ge 0}} A_n$, let $\Proj A$ be the projective 
scheme over $\Spec A_0$, which is associated to $A$. Throughout the paper,
we only consider integral, separated and reduced schemes over $\C$. 
We call them varieties.

Let $X$ be a variety over $\C$.
For a sheaf $\calF$ on $X$ and an open subset $U \subset X$, we denote
the set of local sections of $\calF$ on $U$ by $\calF(U)$ or $\Gamma(U, \calF)$.
We denote the structure sheaf of $X$ by $\calO_X$ and the coordinate ring
of $X$ by $\C[X] = \calO_X(X)$.
For a smooth complex manifold $X$, let $\frD(X)$ be the algebra of 
algebraic differential operators on $X$.

\section{DQ-algebras}
\label{sec:w-algebras}

In this section, we recall the definition of  ($\hbar$-localized)
DQ-algebras according to \cite{KR}.

Let $\hbar$ be an indeterminate. Given $m \in \Z$, let 
$\W_{T^* \C^n} (m)$ be a sheaf of formal series 
$\sum_{k \ge -m} \hbar^k a_k$ ($a_k \in \calO_{T^* \C^n}$)
on the cotangent bundle $T^* \C^n$ of $\C^n$. We set
$\W_{T^* \C^n} = \bigcup_{m} \W_{T^* \C^n} (m)$. We define
a noncommutative $\C((\hbar))$-algebra structure on $\W_{T^* \C^n}$ by
\[
 f \star g = \sum_{\alpha \in \Z_{\ge 0}^n} 
 \hbar^{|\alpha|} \frac{1}{\alpha!} \partial^{\alpha}_{\xi} f
 \cdot \partial^{\alpha}_{x} g
\]
where, for a multi-power $\alpha = (\alpha_1, \dots, \alpha_n) \in \Z^n_{\ge 0}$, 
we set $\alpha! = \alpha_1 ! \cdots \alpha_n !$ and
$|\alpha| = \alpha_1 + \dots + \alpha_n$. Note that $\W_{T^* \C^n}(0)$
is a $\C[[\hbar]]$-subalgebra of $\W_{T^* \C^n}$.

Let $X$ be a complex symplectic manifold with symplectic form $\omega$.
A DQ-algebra on $X$ is a sheaf of $\C((\hbar))$-algebras $\W$ such that
for any point $x \in X$, there is an open neighborhood $U$ of $x$,
a symplectic map $\varphi: U \rightarrow T^* \C^n$, and a $\C((\hbar))$-algebra
isomorphism $\psi: \W \vert_U \isoto \varphi^{-1} \W_{T^* \C^n}$.

We have the following fundamental properties of a DQ-algebra $\W$ as listed in
\cite{KR}.

\begin{enumerate}
\item 
The algebra $\W$ is a coherent and noetherian algebra.
\item
$\W$ contains a canonical $\C[[\hbar]]$-subalgebra $\W(0)$
which is locally isomorphic to $\W_{T^*\C^n}(0)$ (via the maps $\psi$).
We set $\W(m)=\hbar^{-m}\W(0)$.
\item
We have a canonical $\C$-algebra isomorphism
$\W(0)/\W(-1)\isoto \calO_X$ (coming from the canonical isomorphism via the maps
$\psi$).
The corresponding morphism
$\sigma_m: \W(m)\to \hbar^{-m}\calO_X$ is called the {\em symbol map}.
\item We have
\[
\sigma_{0}(\hbar^{-1}[f,g])=\{\sigma_0(f),\sigma_0(g)\}
\]
for any $f$, $g\in\W(0)$. Here $\{\scbul,\scbul\}$
is the Poisson bracket of $X$ which induced from the symplectic structure
of $X$.
\item
The canonical map 
$\W(0)\to \prolim[{m\to\infty}]\W(0)/\W(-m)$
is an isomorphism.
\item
A section $a$ of $\W(0)$ is invertible in $\W(0)$ if and only if
$\sigma_0(a)$ is invertible in $\calO_X$.
\item
Given $\phi$, a $\C((\hbar))$-algebra automorphism of $\W$,
we can find locally an invertible section $a$ of $\W(0)$
such that $\phi=\Ad(a)$.
Moreover $a$ is unique up to a scalar multiple.
In other words, we have canonical isomorphisms
\[
\xymatrix{
\W(0)^\times/\C[[\hbar]]^\times\ar[r]^\sim_{\Ad}\ar[d]_\sim
&\Aut(\W(0))\ar[d]^\sim\\
\W^\times/\C((\hbar))^\times\ar[r]^\sim_{\Ad}&\Aut(\W).
} 
\]
\item
Let $v$ be a $\C((\hbar))$-linear filtration-preserving
derivation of $\W$.
Then there exists locally a  section $a$ of $\W(1)$
such that $v=\ad(a)$.
Moreover $a$ is unique up to a scalar.
In other words, we have an isomorphism
\[
\W(1)/\hbar^{-1}\C[[\hbar]]\isoto[\ad]\Der_{\mathrm{filtered}}(\W).
\]
\item
If $\W$ is a DQ-algebra, then its opposite ring
$\W^\opp$ is a DQ-algebra on $X^\opp$ where $X^\opp$ is
the symplectic manifold with symplectic form $- \omega$.
\end{enumerate}

A tuple $(x_1, \dots, x_n; \xi_1, \dots, \xi_n)$ of elements 
$x_i$, $\xi_i \in \W(0)$ are called quantized symplectic coordinates
of $\W$ if they satisfy $[x_i, x_j] = [\xi_i, \xi_j] = 0$ and
$[\xi_i, x_j] = \hbar \delta_{ij}$. 

Next, we review the notion of F-actions.
Let $X$ be a symplectic manifold with the action of $\bbG_m$:
$\C^* \ni t \mapsto T_t \in \Aut(X)$. We assume there exists
a positive integer $m \in \Z_{> 0}$ such that 
$T^*_t \omega = t^m \omega$ for all $t \in \C^*$.

An F-action with exponent $m$ on $\W$ is an action of 
$\C^*$ on the
$\C$-algebra $\W$,
$\scF_t : T^{-1}_t \W \isoto \W$ for $t \in \C^*$ such that
$\scF_t(\hbar) = t^m \hbar$ and $\scF_t(f)$ depends holomorphically
on $t$ for any $f \in \W$.
An F-action with exponent $m$ on $\W$ extends to an F-action
with exponent $1$ on 
$\W[\hbar^{1/m}] = \C((\hbar^{1/m})) \otimes_{\C((\hbar))} \W$ given by
$\scF_t(\hbar^{1/m}) = t^1 \hbar^{1/m}$. 
For an F-action on $\W$ with exponent $m$, we consider the $\C$-algebra of
global sections of $\W[\hbar^{1/m}]$ which are
 invariant under the F-actions,
and denote it by $\Gamma_{\scF}(\W)$; i.e. 
$\Gamma_{\scF}(\W) \defeq \Gamma(X, \W[\hbar^{1/m}])^{\C^*}$.

For the DQ-algebra $\scW_{T^*\C^n}$ on $T^*\C^n$, we consider an F-action
$\scF_t$ ($t \in \C^*$) defined by $\scF_t(x_i) = t^1 x_i$, 
$\scF_t(\xi_i) = t^1 \xi_i$ and $\scF_t(\hbar) = t^2 \hbar$ for $i=1$, $\dots$,
$n$. It is an F-action of exponent $2$.
 Set $\widetilde{\scW}_{T^*\C^n} = \scW_{T^*\C^n}[\hbar^{1/2}]$. 
Then we have the following injective homomorphism of $\C$-algebras,
\[
 \frD(\C^n) \longrightarrow \widetilde{\scW}_{T^*\C^n}(T^*\C^n), \quad
 x_i \mapsto \hbar^{-1/2} x_i, \quad 
 \frac{\partial}{\partial x_i} \mapsto \hbar^{-1/2} \xi_i.
\]
This homomorphism induces an isomorphism of $\C$-algebras,
$\Gamma_{\scF}(\scW_{T^*\C^n}) \simeq \frD(\C^n)$ 
(see \cite[Lemma 2.9]{BK}).

\section{Quantum Hamiltonian reduction}
\label{sec:setting-q-Hamilton-red}

\subsection{Hamiltonian reduction}
\label{sec:hamilt-reduct}

Let $V$ be a vector space over $\C$. Its cotangent bundle $T^* V$ has
natural symplectic structure. Let $G$ be a reductive algebraic
group which acts algebraically on $V$. 
This action induces a Hamiltonian
 action of $G$ on $T^* V$ and we have
a moment map $\mu_{T^* V} : T^* V \longrightarrow \frg^* = (\Lie G)^*$.
Let $\bbX$ be the set of all characters $G \longrightarrow \C^*$ of $G$
and let $\bbX_\Q = \bbX \otimes_{\Z} \Q$ be the space of fractional
characters.
We fix $\theta \in \bbX_\Q$ and call $\theta$ a stability parameter.
A point $p \in T^* V$ is called $\theta$-semistable if there exist
a function $f \in \C[T^*V]$ and $m \in \Z_{> 0}$ such that
$g \cdot f = \theta(g)^m f$ for any $g \in G$ and $f(p) \ne 0$.
Let $\mu_{\frX} = \mu_{T^* V} \vert_{\frX} : \frX \longrightarrow \frg^*$ be 
the restriction of 
the moment map $\mu_{T^*V}$ to $\frX$.
 The subset $\mu^{-1}_{\frX}(0)$ of $\frX$ is 
closed under the $G$-action.
We consider the affine GIT quotient 
\[
 X_0 = \mu^{-1}_{T^*V}(0) /\!/ G = \Spec \C[\mu_{T^*V}^{-1}(0)]^G,
\]
and the projective GIT quotient 
\[
 X = \mu^{-1}_{T^*V}(0) /\!/_{\theta} G = \Proj \bigoplus_{m \in \Z_{\ge 0}}
 \C[\mu^{-1}_{T^*V}(0)]^{G, \theta^m}.
\]
The inclusion morphism $\frX \longrightarrow T^* V$ induces a morphism
$X \longrightarrow X_0$. 

Throughout this paper, we assume 
\begin{enumerate}
 \item $\mu_{\frX}^{-1}(0)$ is not empty, \label{item:1}
 \item $G$ acts freely on $\mu_{\frX}^{-1}(0)$, \label{item:2}
 \item The moment map $\mu_{T^*V}$ is flat over $\frg^*$, \label{item:3}
 \item The morphism $X \longrightarrow X_0$ is birational and $X_0$ is a normal variety. \label{item:4}
\end{enumerate}
By the assumption \ref{item:1} and \ref{item:2}, the variety $X$ is
isomorphic to the quotient space $\mu_{\frX}^{-1}(0) / G$ and it
is a smooth symplectic manifold. Moreover, we conclude that 
the above morphism $X \longrightarrow X_0$ is a resolution of singularity.
Let $p : \mu^{-1}_{\frX}(0) \longrightarrow X$ be the projection.

Let $\calO_{\frX}$ be the structure sheaf of $\frX$. 
The action of $G$ on $V$ induces an 
equivariant $G$-action
on $\calO_{\frX}$. Moreover, the moment map
$\mu_{\frX}$ induces a comoment map 
$\mu^*_{\frX} : \frg \longrightarrow \calO_{\frX}$. 
Let $\{A_1, \dots, A_{\Dim \frg}\}$ be a basis of $\frg$.
Then we have the isomorphism
$\calO_{\mu_{\frX}^{-1}(0)} \simeq \calO_{\frX} / \sum_{i=1}^{\Dim \frg}\calO_{\frX} \mu^*_{\frX}(A_i)$ 
and the structure sheaf $\calO_{X}$ of $X$ is isomorphic to 
$(p_* (\calO_{\frX} / \sum_{i=1}^{\Dim \frg} \calO_{\frX} \mu^*_{\frX}(A_i)))^G$.

Let $K^{\bullet}(\calO_{\frX}, \{\mu^*_{\frX}(A_1), \dots, \mu^*_{\frX}(A_{\Dim \frg})\})$
be the Koszul complex associated to $\calO_{\frX}$ and the sequence of global
sections $\{\mu^*_{\frX}(A_1), \dots, \mu^*_{\frX}(A_{\Dim \frg})\}$. Then its 
zeroth homology coincides with $\calO_{\mu_{\frX}^{-1}(0)}$.
Note that 
the moment map $\mu_{T^*V}$ is flat  and hence so is $\mu_{\frX}$.
The following lemma is due to M.~Holland (\cite{H}).

\begin{lemma}[\cite{H}, proof of Proposition 2.4]
 \label{lemma:4}
 The sequence of the global sections $\{\mu^*_{\frX}(A_1), \dots, \mu^*_{\frX}(A_{\Dim \frg})\}$
 is a regular sequence in $\calO_{\frX}(\frX)$. Thus, for any open subset
 $\frU$ of $\frX$, higher homology of
 the Koszul complex 
$K^{\bullet}(\calO_{\frX}(\frU), \{\mu^*_{\frX}(A_1), \dots, \mu^*_{\frX}(A_{\Dim \frg})\})$
vanishes; i.e. we have
\[
 H^{n}(K^{\bullet}(\calO_{\frX}(\frU), \{\mu^*_{\frX}(A_1), \dots, \mu^*_{\frX}(A_{\Dim \frg})\})) \simeq
 \begin{cases}
  \calO_{\mu_{\frX}^{-1}(0)}(\frU) & \text{if } n = 0, \\
  0 &  \text{otherwise.}
 \end{cases}
\]
\end{lemma}

\subsection{Quantum Hamilton reduction}
\label{sec:quant-hamilt-reduct}

Let $\frD(V)$ be the ring of algebraic differential operators on $V$.
The action of $G$ on $V$ induces an action on 
$\frD(V)$. Moreover, by differentiating the action of $G$ on $\calO_V$, we have
a morphism of algebras
\[
 \mu_{\frD} : \frg \longrightarrow \frD(V)
\]
which we call a quantized moment map. Fix a parameter $c : \frg \longrightarrow \C$ to be a $G$-invariant linear function, i.e. a character of $\frg$.
We consider the following subquotient
of the algebra $\frD(V)$ with respect to the $G$-action:
\[
 \frD(X_0, c) = (\frL_c)^G, \quad \text{where} \quad
 \frL_c = \frD(V) \bigm/ \sum_{i \in 1}^{\Dim \frg} \frD(V) (\mu_{\frD}(A_i) + c(A_i)).
\]

Let $\W_{T^* V}$ be the standard DQ-algebra associated to the symplectic
structure $T^* V$. Set $\W_{\frX} = \W_{T^* V} \vert_{\frX}$ be its
restriction onto $\frX$. The action of $G$ on $\frD(V)$ induces
an equivariant $G$-action on $\W_{T^* V}$: i.e., for $g \in G$, we
have $T_g \in \Aut(\frX)$ and we have an isomorphism
\[
 \rho_g : T^{-1}_g \W_{\frX} \isoto \W_{\frX}
\]
such that $\rho_g(\hbar) = \hbar$.

The quantized moment map $\mu_{\frD}$ induces the following homomorphism
of algebras
\[
 \mu_{\W} : \frg \xrightarrow{\mu_{\frD}} \frD(V) \hookrightarrow
\W_{T^*V}(T^*V) = \W_{\frX}(\frX).
\]
We call the above homomorphism $\mu_{\W}$ also a quantized moment map.
This homomorphism is known to be satisfied the following properties (see \cite{KR}).
\begin{enumerate}
\item $[\mu_{\W}(A), a] = \frac{d}{dt} \rho_{\exp(tA)}(a) \vert_{t = 0}$,
\item $\sigma_0(\hbar \mu_{\W}(A)) = \mu_{\frX} (A)$,
\item $\mu_{\W}(\textrm{Ad}(g) A) = \rho_g (\mu_{\W}(A))$,
\end{enumerate}
for every $A \in \frg$, $a \in \W_{\frX}$ and $g \in G$.

Let $c$ be a parameter as stated above.
We define the $\W_{\frX}$-module 
$\scL_{c}$ by
\[
 \scL_c = \W_{\frX} \bigm/ \sum_{i=1}^{\Dim \frg}
 \W_{\frX} (\mu_{\W}(A_i) + c(A_i)).
\]
Define a sheaf of algebras $\W_{X,c}$ on $X$ by
\[
 \W_{X,c} = \left(p_{*} \scL_{c} \right)^G.
\]
As shown in \cite{KR}, $\W_{X,c}$ is a DQ-algebra on $X$. 

Set $\scW_{\frX}(m) = \scW_{T^*V}(m) \vert_{\frX}$. Then 
$\{ \scW_{\frX}(m) \}_{m \in \Z}$ is a $\C[[\hbar]]$-algebra filtration of $\scW_{\frX}$
such that $\scW_{\frX}(m) / \scW_{\frX}(m-1) \simeq \hbar^{-m} \calO_{\frX}$ for all $m \in \Z$. 
Define
\[
\scL_c(0) = \W_{\frX}(0) \bigm/
\sum_{i=1}^{r} \W_{\frX}(-1) (\mu_{\W}(A_i) + c(A_i)).
\]
Then $\scL_c(0)$ is a $\scW_{\frX}(0)$-submodule of $\scL_c$ such that
$\scL_c = \scW_{\frX} \otimes_{\scW_{\frX}(0)} \scL_c(0)$. Setting
$\scW_{X,c}(0) = (p_* \scL_c(0))^G$, this sheaf of $\C[[\hbar]]$-algebras
on $X$ gives the $\C[[\hbar]]$-subalgebra in (2) of the list of properties of DQ-algebras
in \refsec{sec:w-algebras}.

Define an F-action on $\scW_{T^* V}$ by $\scF_{t}(x) = t^1 x$, 
$\scF_t(\xi) = t^1 \xi$ and $\scF_t(\hbar) = t^2 \hbar$ where
$x \in V^* \subset \C[T^* V]$, $\xi \in V \subset \C[T^* V]$ and 
$t \in \C^*$. This F-action on $\scW_{T^* V}$ induces an F-action
on $\scW_{X, c}$ with exponent $2$. 
This F-action commutes
with the equivariant $G$-action on $\scW_{T^* V}$; i.e.
$T_g$ and $T_t$ commute,
$\rho_g$ and $\scF_t$ commute, and $\mu_{\scW}(A)$ is $\C^*$-invariant
for $g \in G$, $t \in \C^*$ and $A \in \frg$.

Note that we assume that the moment map $\mu_{T^* V}$ is flat and
$X_0$ has normal singularity.
Then, we have the following proposition.
\begin{proposition}[\cite{BK}, Proposition 3.5]
 \label{prop:3}
 We have the isomorphism of $\C$-algebras 
 $\Gamma_{\scF}(\scW_{X,c}) \simeq \frD(X_0, c)$.
\end{proposition}

\section{Algebras constructed as quantum Hamiltonian reduction}
\label{sec:alg-q-Hamilton}

In this section, we introduce some algebras which are constructed
by quantum Hamiltonian reduction defined in \refsec{sec:setting-q-Hamilton-red}.
Throughout this section, we use the notations introduced in 
\refsec{sec:hamilt-reduct} and \refsec{sec:quant-hamilt-reduct}.

\subsection{Quantization of Quiver varieties}
\label{sec:q-quiver-varieties}

Let $Q = (I, E)$ be a finite quiver. We assume that $Q$ has no loop.
Let $\varepsilon_i \in \Z^I \subset \C^I$ 
be the standard basis corresponding to the vertex $i \in I$. We denote
$\bfv = \sum_{i \in I} v_i \varepsilon_i \in \C^I$ by 
$(v_i)_{i \in I}$.
For $\bfv = (v_i)_{i \in I}$, 
$\bfw = (w_i)_{i \in I} \in \C^I = \Z^I \otimes_{\Z} \C$,
we consider the inner product $\bfv \cdot \bfw = \sum_{i \in I} v_i w_i$.
For $\bfv = (v_i)_{i \in I} \in \C^I$,
set $p(\bfv) = 1 + \sum_{\alpha \in E} v_{out(\alpha)} v_{in(\alpha)} - \bfv \cdot \bfv$.

Fix a dimension vector $\bfv = (v_i)_{i \in I} \in \Z^I$. Let $V$ be a 
vector space $V = \bigoplus_{\alpha \in E} \gHom(\C^{v_{out(\alpha)}}, \C^{v_{in(\alpha)}})$,
and let $T^* V$ be its cotangent bundle. Set $G = \prod_{i \in I} GL(\C^{v_i}) / \C^*_{diag}$,
a reductive algebraic group, 
and let $\frg = \bigoplus_{i \in I} \gl(\C^{v_i}) / \C_{diag}$ be its
Lie algebra where $\C^*_{diag}$ (resp. $\C_{diag}$) is the diagonal subgroup
of $\prod_{i \in I} GL(\C^{v_i})$ (resp. the diagonal Lie subalgebra of
$\bigoplus_{i \in I} \gl(\C^{v_i})$).
 Consider an action of $G$ on $V$ defined by
\[
 g \cdot e = (g_{in(\alpha)} e_\alpha g_{out(\alpha)}^{-1})_{\alpha \in E}
\]
for $g = (g_i)_{i \in I} \in G$ and $e = (e_\alpha)_{\alpha \in E} \in V$.
The action induces a Hamiltonian action on $T^* V$ and a moment map
$\mu_{T^* V} : T^* V \longrightarrow \frg^*$ which is explicitly given by
\begin{gather*}
\mu_{T^* V} : T^* V = \bigoplus_{\alpha \in E} 
 T^* \gHom(\C^{v_{out(\alpha)}}, \C^{v_{in(\alpha)}}) 
 \longrightarrow \frg^* \subset \bigoplus_{i \in I} \gEnd(\C^{v_i}), \\
 \mu_{T^* V}((X^{\alpha}, Y^{\alpha})_{\alpha \in E})
 = \Bigl(\sum_{\substack{\alpha \in E \\ in(\alpha) = i}} X^{\alpha} Y^{\alpha}
 - \sum_{\substack{\alpha \in E \\ out(\alpha) = i}} Y^{\alpha} X^{\alpha}
\Bigr)_{i \in I}.
\end{gather*}
Then we have the following fact.
\begin{lemma}[\cite{C}, Theorem~1.1]
\label{lemma:8}
 The moment map $\mu_{T^*V}$ is a flat morphism if and only if
 $\mu_{T^* V}^{-1}(0)$ has dimension $\bfv \cdot \bfv - 1 + 2 p(\bfv)$.
\end{lemma}

We identify a stability parameter $\theta \in \bbX_\Q$ with
$\theta = (\theta_i)_{i \in I} \in \Q^I$ satisfying $\theta \cdot \bfv = 0$ by
$\theta(g) = \sum_{i \in I} (\det g_i)^{- \theta_i}$ for 
$g = (g_i)_{i \in I} \in G$.
Fix a stability parameter $\theta$ and let $\frX$ be the subset of  
$\theta$-semistable points in $T^* V$ with respect to the $G$-action.
Consider the moment map $\mu_{\frX} = \mu_{T^* V} \vert_{\frX}$
and we have Hamiltonian reductions $X = \mu_{\frX}^{-1}(0) / G$ and
$X_0 = \mu_{T^*V} /\!/ G$
as in \refsec{sec:hamilt-reduct}. By the assumption in \refsec{sec:quant-hamilt-reduct},  $X$ is a smooth symplectic manifold.

Next we consider quantum Hamiltonian reduction. The action of $G$
on $V$ induces an action on $\frD(V)$ and $\scW_{\frX}$.
By differentiating the
action of $G$ on $\calO_V$, we have quantized moment maps
$\mu_{\frD} : \frg \longrightarrow \frD(V)$ and
$\mu_{\scW} : \frg \longrightarrow \scW_{\frX}(\frX)$. 
As written in \cite[Lemma 3.1]{H}, the quantized moment maps are given explicitly by
\begin{equation}
 \mu_{\scW}(A^{(i)}_{pq}) = 
  \sum_{\substack{\alpha \in E \\ out(\alpha) = i}} \sum_{j=1}^{v_{in(\alpha)}}
  \hbar^{-1} x^{\alpha}_{jp} \xi^{\alpha}_{jq} - 
\sum_{\substack{\alpha \in E \\ in(\alpha) = i}} \sum_{j=1}^{v_{out(\alpha)}}
  \hbar^{-1} x^{\alpha}_{qj} \xi^{\alpha}_{pj} 
\end{equation}
where $A^{(i)}_{pq}$ is the $(p,q)$-th matrix unit of the $i$-th summand $\gl(\C^{v_i})$ of $G$ and
$(x^{\alpha}_{pq}; \xi^{\alpha}_{pq})_{\substack{1 \le p \le v_{in(\alpha)}, \\ 1 \le q \le v_{out(\alpha)}}}$ is the standard symplectic coordinates of
the symplectic vector space $T^* \gHom(\C^{v_{out(\alpha)}}, \C^{v_{in(\alpha)}})$. 
Fix a parameter 
$c : \frg \longrightarrow \C$ as in \refsec{sec:quant-hamilt-reduct}. 
It is easy to see that the parameter $c$ can be identified with
$c = (c_i)_{i \in I} \in \C^I$ satisfying $c \cdot \bfv = 0$
by $c = \sum_{i \in I} c_i \Tr_{\gl(\C^{v_i})}$.
Then we have quantum Hamiltonian reductions of $\frD(V)$ and
$\scW_{\frX}$ with respect to the $G$-action as follows:
\begin{align*}
 \frD(X_0, c) &= 
\bigl(\frD(V) \bigm/ \sum_{i \in I} \sum_{p,q = 1}^{v_i}
 \frD(V)( \mu_{\frD}(A^{(i)}_{pq}) + c_i \delta_{pq}) \bigr)^G, \\
 \scW_{X,c} &= \bigl(\scW_{\frX} \bigm/ \sum_{i \in I} \sum_{p,q = 1}^{v_i}
 \scW_{\frX}( \mu_{\scW}(A^{(i)}_{pq}) + c_i \delta_{pq}) \bigr)^G.
\end{align*}

\subsubsection{Deformed preprojective algebras}
\label{sec:deform-preproj-alg}

Let $Q = (I, E)$ be a quiver whose underlying diagram $Q_0$ is a Dynkin 
diagram of type affine ADE. We identify the lattice $\Z^I$ with the
root lattice associated to the Dynkin diagram $Q_0$. We also
identify $\varepsilon_i \in \Z^I$
with the simple root associated to the vertex $i \in I$. 

Set the dimension
vector $\bfv = \delta$ where $\delta$ is the minimal positive imaginary 
root of $Q_0$. Fix a stability parameter $\theta \in \Q^I$
such that $\theta \cdot \delta = 0$ and 
$\theta \cdot \alpha \ne 0$ for any
root $\alpha$ satisfying $\alpha_i \le 1$ for all
$i \in I$ and $\alpha \ne \delta$.
Applying the facts in \refsec{sec:q-quiver-varieties} to the above $Q$,
$\bfv$ and $\theta$, we have the Hamiltonian reduction $X$ and
the quantum Hamiltonian reduction $\scW_{X,c}$ on $X$ for a parameter
$c \in \C^I$ such that $c \cdot \delta = 0$.
Then we have the following facts.

\begin{lemma}[\cite{C}, \cite{H}]
\label{lemma:9}
 The moment map 
 $\mu_{T^* V} : T^*V \longrightarrow \frg^*$ is a flat morphism.
\begin{proof}
 It is easy to check the equality
 $\Dim \mu^{-1}_{T^*V}(0) = \delta \cdot \delta - 1 + 2 p(\delta)$.
 Therefore, by \reflemma{lemma:8}, $\mu_{T^*V}$ is a flat morphism.
\end{proof}
\end{lemma}

\begin{proposition}[\cite{K}, Corolary 3.12]
 The Hamiltonian reduction $X$ is a smooth symplectic manifold and it is a 
 minimal resolution of the Kleinian singularity of type $Q_0$.
\end{proposition}

The quantum Hamiltonian reduction $\frD(X_0, c)$ was studied by Holland
in \cite{H}. He showed the algebra $\frD(X_0, c)$ is isomorphic to the 
deformed preprojective algebra corresponding to the quiver $Q$ which was
introduced by Crawley-Boevey and Holland in \cite{CBH}. 
On the other hand, we have the isomorphism of \refprop{prop:3}. Then,
we have the following proposition.

\begin{proposition}[\cite{H}]
 The algebra $\Gamma_{\scF}(\scW_{X,c})$ is isomorphic to the 
 deformed preprojective algebra $\calO^{c + \varepsilon_0 + \partial}$
 introduced by Crawley-Boevey and Holland in \cite{CBH},
 where $\partial = (\partial_i)_{i \in I} \in \Z^I$ is defined by
 $\partial_i = - \delta_i + \sum_{\substack{\alpha \in E \\ out(\alpha) = i}} \delta_{in(\alpha)}$.
 \begin{proof}
  By \refprop{prop:3}, we have an isomorphism 
  $\Gamma_{\scF}(\scW_{X,c}) \simeq \frD(X_0, c)$. \cite[Corollary 4.7]{H}
  says that there exists an isomorphism between
  $\frD(X_0, c)$ and the deformed preprojective algebra 
  $\calO^{c + \varepsilon_0 + \partial}$ and hence we have the isomorphism of
  the proposition.
 \end{proof}
\end{proposition}

\subsubsection{Symplectic reflection algebras}
\label{sec:RCA}

Let $Q' = (I', E')$ be a quiver whose underlying diagram $Q'_0$ is a
Dynkin diagram of type affine ADE. We identify $\Z^{I'}$ with the root
lattice associated to the Dynkin diagram $Q'_0$ and $\varepsilon_i \in \Z^{I'}$
with the simple root corresponding to the vertex $i \in I'$.
Let $Q = (I = I' \sqcup \{\infty\} , E = E' \sqcup \{\infty \rightarrow 0\})$
be a quiver defined by adding a vertex $\infty$ and an arrow 
$\infty \rightarrow 0$ where $0$ is the extended vertex of $Q'_0$. The quiver
$Q$ is called a Calogero-Moser quiver.

Set a dimension vector $\bfv = n \delta + \varepsilon_{\infty}$ where
$\delta$ is the minimal positive imaginary root of $Q'_0$. 
Fix a stability parameter $\theta = (\theta_i)_{i \in I'} \in \Q^{I'}$ 
such that $\theta \cdot \alpha \ne 0$ for any root 
$\alpha = (\alpha_i)_{i \in I'}$ satisfies $\alpha_i \le n$ for all $i \in I'$.
Then we have the Hamiltonian reduction $X$ and the quantum Hamiltonian 
reductions $\frD(X_0, c)$ and $\scW_{X,c}$ for a parameter $c = (c_i)_{i \in I'} \in \C^{I'}$.
By combining \cite[Theorem 3.2.3(iii)]{GG2} and \reflemma{lemma:8}, we have
the following lemma (see also \cite[Theorem 2.6]{Go2}).

\begin{lemma}[\cite{GG2}, \cite{C}]
 The moment map $\mu_{T^*V}$ is a flat morphism.
\end{lemma}

\begin{proposition}[\cite{Na} Theorem 2.8. See also Theorem 4.1]
 The Hamiltonian reduction $X$ is a smooth symplectic manifold of dimension
 $2 n$. Moreover, we have a resolution of singularity 
 $X \longrightarrow \C^{2n} / W$ where $W$ is the wreath product 
 $W = \frS_n \wr \Gamma$ of the finite subgroup $\Gamma$ of $SL_2(\C)$
 corresponding to the Dynkin diagram $Q'_0$.
\end{proposition}

The quantum Hamiltonian reduction $\frD(X_0, c)$ was studied in
\cite{GG2}, \cite{Go2} and \cite{EGGO}, and is known to be isomorphic to
a symplectic reflection algebra of corresponding type, which was
introduced by Etingof and Ginzburg in \cite{EG}. 
Together with
the isomorphism of \refprop{prop:3}, we conclude the following proposition.

\begin{proposition}[\cite{GG}, \cite{Go2}, \cite{EGGO}]
 The algebra $\Gamma_{\scF}(\scW_{X, c})$ is isomorphic to a symplectic
 reflection algebra corresponding to the quiver $Q$.
\end{proposition}

\subsection{Quantization of hypertoric varieties}
\label{sec:q-hypertoric-var}

In this subsection, we review hypertoric varieties and their quantization
introduced by Musson and Van den Bergh in \cite{MVdB}. For facts about
hypertoric varieties, see also \cite{Proudfoot}.

Let $V = \C^n$ be an $n$-dimensional vector space and $G = (\C^*)^d$
be a $d$-dimensional torus acting algebraically on $V$. We take a basis
$\{v_1, \dots, v_n\}$ such that there exists a matrix 
$M = (\mu_{ij})_{\substack{1 \le i \le d \\ 1 \le j \le n}}$ and $G$ acts by
\[
 (t_1, \dots, t_d) \cdot v_i = t_1^{\mu_{1i}} \dots t_d^{\mu_{di}} v_i
\]
for $(t_1, \dots, t_d) \in G$. We assume that $M$ is a unimodular matrix.
Consider the cotangent bundle $T^* V$ of $V$. We denote the coordinates of
$V$ with respect to the basis $\{v_1, \dots, v_n\}$ by 
$(x_1, \dots, x_n)$; i.e. $x_i(v_j) = \delta_{ij}$. Let $(\xi_1, \dots, \xi_n)$
be the dual coordinates so that the symplectic form on $T^* V$
is given by $\omega = \sum_{i = 1}^n d x_i \wedge d \xi_i$. The action of
$G$ on $V$ induces a Hamiltonian action on $T^* V$ and we have a moment
map $\mu_{T^* V}$ with respect to this action.

\begin{lemma}[\cite{BK}, Lemma~4.7]
 The moment map $\mu_{T^* V}$ is flat.
\end{lemma}

We identify a stability parameter $\theta \in \bbX_\Q$ with
$(\theta_1, \dots, \theta_d) \in \Q^d$ and 
fix it. Let $\frX$ be the set of 
$\theta$-semistable points of $T^* V$. Set $\mu_{\frX} = \mu_{T^* V} \vert_{\frX}$. 
Then we have Hamiltonian reductions of these spaces:
\[
 X = \mu^{-1}_{T^*V}(0) /\!/_{\theta} G = \mu^{-1}_{\frX}(0) / G, \qquad \text{and} \qquad
 X_0 = \mu^{-1}_{T^* V}(0) /\!/ G.
\]
For $j = 1$, $\dots$, $n$, set $\mu_j = {}^t(\mu_{1j}, \dots, \mu_{dj}) \in \C^d$. Note that we assume $M = (\mu_{ij})_{\substack{1 \le i \le d \\ 1 \le j \le n}}$ is a unimodular matrix. Then we have the following proposition.
\begin{proposition}[\cite{BK}, Corollary 4.13]
 If the stability parameter $\theta \in \bbX_\Q$
 satisfies $\theta \not\in \sum_{j \in J} \Q \mu_j$ for all
 $J \subset \{1, \dots, n\}$ such that 
 $\Dim_{\Q}(\Span_{\Q}\{\mu_j \}_{j \in J}) = d - 1$, 
 then the variety $X$ is a smooth symplectic manifold and
 we have the resolution of singularity $X \longrightarrow X_0$.
\end{proposition}

Consider the algebra of algebraic differential operators $\frD(V)$ on $V$.
The action of $G$ on $V$ induces an action on $\frD(V)$. Moreover,
by differentiating the action of $G$ on $\calO_V$, we have a quantized moment map $\mu_{\frD}$ 
explicitly as follows:
\[
 \mu_{\frD} : \frg = \bigoplus_{i=1}^d \C A_i \longrightarrow \frD(V), \qquad
 \mu_{\frD}(A_i) = \sum_{j=1}^n \mu_{ij} x_j \frac{\partial}{\partial x_j}.
\]
We identify a character $c : \frg \longrightarrow \C$ of $\frg$ with 
$(c_1, \dots, c_d) \in \C^d$ and fix it.
Set 
\[
 \frL_c = \frD(V) \bigm/ \sum_{i=1}^d (\mu_{\frD}(A_i) + c_i),
\]
and we have a quantum Hamiltonian reduction $\frD(X_0, c) = (\frL_c)^G$ of $\frD(V)$ with
respect to the $G$-action. 
This algebra $\frD(X_0,c)$ is known as an algebra
studied by Musson and Van den Bergh in \cite{MVdB}.

Next we consider the canonical DQ-algebra $\scW_{T^*V}$ on $T^*V$.
The action of $G$ on $V$ induces an equivariant action on $\scW_{T^*V}$
and on $\scW_{\frX} = \scW_{T^*V} \vert_{\frX}$. The quantized moment
map $\mu_{\frD}$ induces the following quantized moment map for
$\scW_{T^* V}$,
\begin{align*}
 \mu_{\scW}: \frg &\longrightarrow \frD(V) \hookrightarrow \scW_{\frX}(\frX), \\
 \mu_{\scW}(A_i) &= \sum_{j=1}^n \hbar^{-1} \mu_{ij} x_j \xi_j.
\end{align*}
For the above parameter $c = (c_1, \dots, c_d) \in \C^d$, set
\[
 \scL_c = \scW_{\frX} \bigm/ \sum_{i=1}^d \scW_{\frX} (\mu_{\scW}(A_i) + c_i).
\]
Then we have a quantum Hamiltonian reduction
$\scW_{X, c} = (p_* \scL_c)^G$.
The following proposition is a direct consequence of \refprop{prop:3}.

\begin{proposition}[\cite{BK}, Proposition 3.5]
 We have an isomorphism of algebras
 $\Gamma_{\scF}(\scW_{X,c}) \simeq \frD(X_0, c)$
and this algebra is the algebra studied by Musson and Van den Bergh
in \cite{MVdB}.
\end{proposition}

\section{BRST cohomology}
\label{sec:brst-cohomology}

In this section, we review the definition of BRST cohomologies
in terms of graded superalgebras. First we recall notions
of superalgebras and Clifford algebras.

\subsection{Superalgebras and Clifford algebras}
\label{sec:super-algebr-cliff}

A superalgebra $R$ is a $\Z/2\Z$-graded algebra over $\C$,
$R = R_0 \oplus R_1$. 
A homogeneous element in $R_0$ is called
an even element, one in $R_1$ is called an odd element. 
For a homogeneous element $a$, we denote the $\Z/2\Z$-degree of $a$
by $|a|$ and call it the ``parity'' of $a$, instead of ``degree''.
For homogeneous elements $a$, $b \in R$, we define 
super-commutator of $a$ and $b$ 
by $[a, b] = a b - (-1)^{|a| |b|} b a$ and extend it
onto whole $R$ linearly.

In this paper, we consider $\Z$-graded superalgebras. A $\Z$-graded 
superalgebra is a superalgebra with (usual) $\Z$-grading. 
In the rest of paper, we always use the notion ``degree'' for the
$\Z$-grading and call the $\Z/2\Z$-grading ``parity''.

For two superalgebras $R$ and $S$, the tensor product of these
algebras is a superalgebra $R \otimes_{\C} S$ with products
$(a \otimes b) \cdot (a' \otimes b') = (-1)^{|a'| |b|} a a' \otimes b b'$.
Its parity is given by
\[
 (R \otimes S)_0 = R_0 \otimes S_0 \oplus R_1 \otimes S_1, \quad
 (R \otimes S)_1 = R_1 \otimes S_0 \oplus R_0 \otimes S_1.
\]

For a finite-dimensional
 vector space $\frg$ with a basis $\{A_1, \dots, A_{\Dim \frg}\}$, 
a Clifford algebra $Cl(\frg \oplus \frg^*)$ associated to the 
symplectic vector space $\frg \oplus \frg^*$ is a superalgebra generated by
$2 \Dim \frg$ odd elements $\psi_1$, $\dots$, $\psi_{\Dim \frg}$, $\psi^*_1$, $\dots$,
$\psi^*_{\Dim \frg}$ with the following defining relations
\[
 [\psi_i, \psi_j] = [\psi^*_i, \psi^*_j] = 0, \quad
 [\psi_i, \psi^*_j] = \delta_{i,j} 1 \quad \text{for } 1 \le i, j \le \Dim \frg.
\]
We regard the exterior algebra $\Lambda(\frg)$ (resp. $\Lambda(\frg^*)$) as
a subalgebra of $Cl(\frg \oplus \frg^*)$ by $A_i \mapsto \psi_i$ 
(resp. $A^*_i \mapsto \psi^*_i$ where $\{A^*_1, \dots, A^*_{\Dim \frg}\}$ 
is the dual basis).
As a vector space, the Clifford algebra $Cl(\frg \oplus \frg^*)$ is isomorphic to the
tensor product of these exterior algebras, i.e. 
$Cl(\frg \oplus \frg^*) \simeq_{\C} \Lambda(\frg) \otimes_{\C} \Lambda(\frg^*)$.
The Clifford algebra is $\Z$-graded by degree $\deg(\psi_i) = -1$ 
and $\deg(\psi^*_i) = 1$ for $i=1$, $\dots$, $\Dim \frg$. We have
\begin{gather*}
 Cl(\frg \oplus \frg^*) = \bigoplus_{n \in \Z} Cl^n(\frg \oplus \frg^*), \\
 Cl^n(\frg \oplus \frg^*) = \begin{cases}
	     \bigoplus_{i+j=n} \Lambda^{-i}(\frg) \otimes_{\C} \Lambda^j(\frg^*) &
	     \text{for } - \Dim \frg \le n \le \Dim \frg, \\
	     0 & \text{otherwise.}
	    \end{cases}
\end{gather*}
Moreover $Cl(\frg \oplus \frg^*)$ is a double-graded superalgebra as follows:
\begin{gather*}
 Cl(\frg \oplus \frg^*) = \bigoplus_{i,j \in \Z} Cl^{i,j}(\frg \oplus \frg^*), \\
 Cl^{i,j}(\frg \oplus \frg^*) = \begin{cases}
	     \Lambda^{-i}(\frg) \otimes_{\C} \Lambda^j(\frg^*) &
	     \text{for } - \Dim \frg \le i \le 0, 0 \le j \le \Dim \frg, \\
	     0 & \text{otherwise.}
	    \end{cases}
\end{gather*}

\subsection{Definition of BRST cohomologies}
\label{sec:def-BRST-alg}

In this subsection, we define the BRST cohomologies of $\scW_{\frX}$ and $\frD(V)$
with respect to the $G$-action. In this subsection, we use the notation
introduced in \refsec{sec:quant-hamilt-reduct}.

Fix $\frU$ be an open subset of $\frX$ which is closed under the $G$-action.
Set $R = \scW_{\frX}(\frU)$ or $\frD(V)$. 
Consider the Clifford
algebra $Cl(\frg \oplus \frg^*)$. We regard $R$ as
a superalgebras with purely even parity; i.e. $R_0 = R$ and 
$R_1 = 0$.
Consider the tensor product superalgebra
$C(R) = R \otimes_{\C} Cl(\frg \oplus \frg^*)$ of 
$R$ and $Cl(\frg \oplus \frg^*)$.
The $\Z$-grading of $Cl(\frg \oplus \frg^*)$ induces a $\Z$-grading of $C(R)$ as follows;
\[
 C(R) = \bigoplus_{n \in \Z} C^n(R), \quad
 C^n(R) = R \otimes_{\C} Cl^n(\frg \oplus \frg^*).
\]
For a character $c: \frg \longrightarrow \C$, 
let $Q_c$ be an odd element of $C(R)$ defined by 
\[
 Q_c = \sum_{i = 1}^{\Dim \frg} (\mu_{\scW}(A_i) + c(A_i)) \otimes \psi^*_i
 + \frac{1}{2} \sum_{i,j,k} \chi_{ij}^k \otimes \psi_k \psi^*_i \psi^*_j
\]
where $\chi_{ij}^{k}$ is the structure constant
of $\frg$; namely $[A_i, A_j] = \sum_k \chi_{ij}^k A_k$.
Note that the F-action on $\scW_{\frX}$ induces a $\C^*$-action on $C(R)$ and
$Q_c$ is a $\C^*$-invariant element with respect to this action.
The adjoint operator 
$\ad Q_c = [Q_c, \bullet\,]$ on $C(R)$ is 
a homogeneous operator with degree $+ 1$. Since it is an
odd operator, we have $(\ad Q_c)^2 = 0$. Thus,
the pair $(C(R), \ad Q_c)$ is a cochain complex.
We call it a BRST complex of $R$ with respect to the action of $\frg$ 
and its cohomologies
\[
 H_{BRST, c}^{\bullet}(\frg, R) = H^{\bullet}(C(R), \ad Q_c)
\]
are called BRST cohomologies of $R$ with respect to the action of $\frg$. 
Note that the algebra structure
of $R$ induces a graded algebra structure on $H_{BRST, c}^{\bullet}(\frg, R)$.

Next, we introduce BRST cohomology of a sheaf of algebras.
For an open subset $U$ of $X$, let $\frU$ be an open subset of $\frX$ such that
$\frU$ is closed under the $G$-action and
$p^{-1}(U) = \frU \cap \mu_{\frX}^{-1}(0)$.
Then a BRST cohomology sheaf $\calH^{\bullet}_{BRST, c}(\frg, \scW_{\frX})$ of 
$\scW_{\frX}$ with respect to $\frg$ is defined as 
a sheaf on $X$ which is obtained by pushing-forward of the sheaf
 associated
to the presheaf $\frU \mapsto H^{\bullet}_{BRST, c}(\frg, \scW_{\frX}(\frU))$.

\begin{lemma}[cf. \cite{AKM}, Theorem 1.3.2.1]
\label{lemma:7}
 Assume that the moment map $\mu_{T^*V}$ is flat.
 The sheaf $\frU \mapsto H^{\bullet}_{BRST, c}(\frg, \scW_{\frX}(\frU))$
 is supported on $\mu_{\frX}^{-1}(0)$ and hence 
 $\calH^{\bullet}_{BRST,c}(\frg, \scW_{\frX})(U)$ does not depend on
 the choice of $\frU$. Namely, 
 $\calH^{\bullet}_{BRST,c}(\frg, \scW_{\frX})$ is a well-defined sheaf on
 $X$.
\end{lemma}

This lemma will be proved in \refsec{sec:vanish-neg-cohom}.

\section{Computing BRST cohomologies}
\label{sec:comp-brst-cohom}

Throughout this section, we use the notations introduced in
\refsec{sec:setting-q-Hamilton-red} and \refsec{sec:brst-cohomology}.
Fix an open subset $U$ of $X$ and an open subset $\frU$ of $\frX$ such that
$\frU$ is closed under the $G$-action and
$p^{-1}(U) = \frU \cap \mu^{-1}_{\frX}(0)$. Set $R = \scW_{\frX}(\frU)$.

\subsection{Double complex}
\label{sec:double-complex}

We consider a double complex structure of 
the BRST complex $(C(R) = R \otimes Cl(\frg \oplus \frg^*), \ad Q_c)$.
Set
\[
 C(R) = \bigoplus_{m,n \in \Z} C^{m,n}(R), \qquad
 C^{m,n}(R) = R \otimes Cl^{m,n}(\frg \oplus \frg^*),
\]
and
\begin{align}
\label{eq:1}
 d_+ : C^{m,n}(R) &\longrightarrow C^{m,n+1}(R), \\
 d_+ (a \otimes \varphi \varphi^*) &=
 \sum_{i=1}^{\Dim \frg} [\mu_{\scW}(A_i), a] \otimes \varphi \psi_i^* \varphi^* 
 + \frac{1}{2} \sum_{i,j,k} \chi_{ij}^k a \otimes \varphi [\psi_k, \varphi^*] 
 \psi_i^* \psi_j^* \nonumber\\
 d_- : C^{m,n}(R) &\longrightarrow C^{m+1,n}(R), \nonumber\\
 d_- (a \otimes \varphi \varphi^*) &=
 \sum_{i=1}^{\Dim \frg} a \star (\mu_{\scW}(A_i) + c(A_i)) \otimes [\psi_i^*, \varphi] \varphi^*
 + \frac{1}{2} \sum_{i,j,k} \chi_{ij}^{k} a \otimes \psi_k [\psi_i^* \psi_j^*, \varphi] \varphi^* \nonumber
\end{align}
where $a \in R$, $\varphi \in \Lambda^{-m}(\frg)$
and $\varphi^* \in \Lambda^n(\frg^*)$.
By direct calculation, $(C(R), d_+, d_-)$ is a double
complex and $\ad Q_c = d_+ + d_-$. Moreover, by considering a spectral 
sequence associated to the double complex $(C(R), d_+, d_-)$, we have the
following lemma.

\begin{lemma}
\label{lemma:3}
 Consider a spectral sequence $E_r^{p,q}$ associated to the double complex 
 $(C(R), d_+, d_-)$ whose second term is given by 
 $E_2^{p,q} = H^p(H^q(C(R), d_-), d_+)$. Then
 the spectral sequence converges to the total cohomology
 $H^{p+q}(C(R), \ad Q_c) = H^{p+q}_{BRST, c}(\frg, R)$.
\begin{proof}
 Since the double complex $C(R) = \bigoplus_{m,n} C^{m,n}(R)$ is bounded
 on $- \Dim \frg \le m \le 0$ and $0 \le n \le \Dim \frg$, the spectral
 sequence is convergent and weekly convergent in the sense of \cite[Chapter 15]{CE}.
\end{proof}
\end{lemma}

\subsection{Vanishing of negative BRST cohomologies}
\label{sec:vanish-neg-cohom}

In this section, we show that the negative degrees of the BRST
cohomology $\calH^{\bullet}_{BRST, c}(\frg, \scW_{\frX})$ vanish.

Consider the $\C[[\hbar]]$-algebra filtration $\{\scW_{\frX}(m)\}_{m \in \Z}$ of $\scW_{\frX}$.
The following properties of the filtration immediately follow from its definition.

\begin{lemma}[\cite{BK}, Lemma~2.2]
\label{lemma:1}
 The filtration $\{\scW_{\frX}(m)\}_{m \in \Z}$ is exhaustive and separated;
 i.e.
\begin{enumerate}
 \item exhaustive: $\bigcup_{m \in \Z} \scW_{\frX}(m) = \scW_{\frX}$, and
 \item separated: $\bigcap_{m \in \Z} \scW_{\frX}(m) = 0$.
\end{enumerate}
\end{lemma}

We introduce a filtration $\{F_m C(R)\}_{m \in \Z}$ of $C(R)$ as follows:
For $m \in \Z$ and $i$, $j \in \Z$, set
\[
 F_{m} C^{i,j}(R) = \scW_{\frX}(m+i)(\frU)
 \otimes_{\C} \Lambda^{-i}(\frg) \otimes_{\C} \Lambda^j(\frg^*).
\]
It is easy to check that we have 
$d_-(F_m C^{i,j}(R)) \subset F_m C^{i+1,j}(R)$ and hence 
the filtration is a filtration of complex $(C^{\bullet, j}(R), d_-)$
for each $j \in \Z$. The following lemma immediately follows from
\reflemma{lemma:1}.

\begin{lemma}
\label{lemma:2}
 The filtration $\{F_m C^{\bullet, j}(R)\}_{m \in \Z}$ is exhaustive
 and separated for each $j \in \Z$.
\end{lemma}

Fix $j \in \Z$ such that $0 \le j \le \Dim \frg$.
The filtration $\{F_m C^{\bullet, j}(R)\}_{m \in \Z}$ induces 
a spectral sequence $\pE^r_{p,q}$ as follows:
\begin{align*}
 \pZ_{p,q}^r &= F_p C^{p+q,j}(R) \cap d_{-}^{-1}(F_{p-r} C^{p+q+1,j}(R)), \\
 \pB_{p,q}^r &= F_p C^{p+q,j}(R) \cap d_{-}(F_{p+r} C^{p+q-1,j}(R)), \\
 \pE_{p,q}^r &= \pZ_{p,q}^r / \pB_{p,q}^r + \pZ_{p+1,q-1}^{r-1},
\end{align*}
for $r \in \Z_{\ge 0}$, $p$, $q \in \Z$. The differential
$d_{-}^{(r)}: \pE_{p,q}^r \longrightarrow \pE_{p+r, q-r+1}^r$ is naturally
induced from the original differential $d_{-}$.

By \reflemma{lemma:2}, we have the following lemma.
\begin{lemma}
 \label{lemma:6}
 The spectral sequence $\pE_{p,q}^r$ is convergent and weekly convergent
 to the cohomology
 $\gr_p H^{p+q}(C^{\bullet,j}(R), d_{-})$ for each $j \in \Z$ in the sense
 of \cite[Chapter 15]{CE}.
\begin{proof}
 Since the filtration is exhaustive and separated (\reflemma{lemma:2}),
 the spectral sequence $\pE^r_{p,q}$ is weekly convergent.
 By \reflemma{lemma:2}, we have
 $\bigcap_p F_p C^{n,j} = 0$ for any $j \in \Z$. By definition 
 $F_p H^n(C^{\bullet,j}(R)) = (F_p C^{n,j}(R) \cap \Ker d_{-} + \Im d_{-}) / \Im d_{-}$, and we have
 $\bigcap_{p} F_p H^n(C^{\bullet,j}(R)) = 0$. Hence, the spectral sequence
 $\pE^r_{p,q}$ is convergent. 
\end{proof}
\end{lemma}

We have $\pE_{p,q}^0 \simeq \gr_p C^{p+q,j}(R)$ for each $p$, $q \in \Z$.
Moreover, the complex 
$(\pE_{p,q}^0, d_{-}^{(0)}) \simeq (\gr_p C^{p+q,j}(R), d_-)$ is isomorphic to
the Koszul complex \\
$K^{- (p+q)}(\calO_{\frX}(\frU), \{\mu_{\frX}^*(A_1), \dots, \mu_{\frX}^*(A_{\Dim \frg})\})$. 
Thus, by \reflemma{lemma:4}, we have the following facts.

\begin{lemma}
 \label{lemma:5}
 For each $j \in \Z$ and $p$, $q \in \Z$, we have an isomorphism
 \[
  \pE_{p,q}^1 \simeq H^{p+q}(\gr_p C^{\bullet, j}(R), d_{-})
 \simeq
 \begin{cases}
  \hbar^{-p} \calO_{\mu_{\frX}^{-1}(0)}(\frU) \otimes_{\C} \Lambda^j(\frg) 
  & \text{if } p+q = 0, \\
  0 & \text{otherwise.}
 \end{cases}
 \]
\end{lemma}

\begin{proposition}
\label{prop:1}
 The spectral sequence $\pE_{p,q}^r$ collapses at $r=1$ and we have
\[
 H^{n}(C^{\bullet, j}(R), d_{-}) \simeq 
\begin{cases}
 \scL_c(\frU) \otimes_{\C} \Lambda^j(\frg^*) & \text{if } n = 0, \\
 0 & \text{otherwise.}
\end{cases}
\]
\begin{proof}
 By \reflemma{lemma:5}, we have $\pE_{p,q}^1 = 0$ for $p+q \ne 0$ and hence
 the differential 
 $d_{-}^{(1)} : \pE_{p,q}^1 \longrightarrow \pE_{p+1,q}^1$ is zero morphism
 for all $p$, $q \in \Z$. Thus, the spectral sequence $\pE_{p,q}^{r}$ collapses
 at $r=1$. Then, by \reflemma{lemma:6}, we have
\begin{multline*}
 H^{p+q}(C^{\bullet, j}(R), d_{-}) \simeq \bigoplus_{p \in \Z}
 \pE_{p,q}^{\infty} \simeq
 \bigoplus_{p \in \Z}\pE_{p,q}^1 \simeq \\
 \bigoplus_{p \in \Z} \hbar^{-p} \calO_{\mu^{-1}_{\frX}(0)}(\frU) \otimes_{\C} \Lambda^{j}(\frg^*) \simeq 
 \scL_c(\frU) \otimes_{\C} \Lambda^j(\frg^*).
\end{multline*}
\end{proof}
\end{proposition}

{\em Proof of \reflemma{lemma:7}}: 
By \refprop{prop:1}, the spectral sequence $E^{p,q}_r$ collapses at $r=2$ and 
we have the isomorphism
\[
 H_{BRST,c}^n(\frg, \scW_{\frX}(\frU)) \simeq H^n(\bigoplus_j H^0(C(R), d_-), d_+) 
 \simeq H^n(\scL_c(\frU) \otimes \Lambda^{\bullet}(\frg), d_+).
\]
By the property (6) in \refsec{sec:w-algebras},
the sheaf $\scL_c$ is supported on $\mu_{\frX}^{-1}(0)$ and we have
the claim of \reflemma{lemma:7}.

\subsection{Positive BRST cohomologies}
\label{sec:pos-BRST-cohom}

By \reflemma{lemma:7}, the sheaf $\scL_c$ with the equivariant $G$-action
is supported on $\mu^{-1}_{\frX}(0)$. Since 
$p^{-1}(U) = \mu_{\frX}^{-1}(0) \cap \frU$ is closed under the $G$-action,
$\scL_c(\frU)$ is a $G$-module.
The $\scW_{\frX}(\frU)$-module $\scL_c(\frU)$ has a filtration as
a $\scW_{\frX}(0)(\frU)$-module $\{ \scL_c(m)(\frU) \}_{m \in \Z}$.
This filtration coincides with the filtration of the cohomology
$H^{0}(C^{\bullet, j}, d_{-}) \simeq \scL_c(\frU) \otimes_{\C} \Lambda^j(\frg^*)$; Namely, we have
$F_m H^{0}(C^{\bullet, j}, d_{-}) \simeq \scL_c(m)(\frU) \otimes_{\C} \Lambda^j(\frg^*)$ for each $0 \le j \le \Dim \frg$. 
Since the $G$-module $\scL_c(\frU)$ over $\C$ is completely reducible and
the filtration is a filtration as a $G$-module,
we have the isomorphism $\scL_c(\frU) \simeq \gr \scL_c(\frU)$ as 
$G$-modules. Moreover, by \refprop{prop:1}, we have
\[
 \scL_c(\frU) \simeq \gr \scL_c(\frU) \simeq \calO_{\mu^{-1}_{\frX}(0)}(\frU) \otimes_{\C} \C((\hbar))
\]
as $G$-modules, where $G$ acts trivially on $\C((\hbar))$.
Now we assume that the $G$-bundle 
$\mu^{-1}_{\frX}(0) \longrightarrow X$ is (locally) trivial on the open
subset $\frU$. Then
we have $\calO_{\mu^{-1}_{\frX}(0)}(\frU) \simeq \calO_{\mu^{-1}_{\frX}(0)}(\frU)^G \otimes_{\C} \C[G]$
as $G$-modules, where $G$ acts trivially on $\calO_{\mu^{-1}_{\frX}(0)}(\frU)^G$ and
by left translation on $\C[G]$. Therefore, we conclude that 
\begin{equation}
\label{eq:2}
\scL_c(\frU) \simeq \scL_c(\frU)^G \otimes_{\C} \C[G]
\end{equation}
as $G$-modules, where $G$ acts trivially on $\scL_c(\frU)^G$.

\begin{proposition}
\label{prop:2}
 For an open subset $\frU$ on which the $G$-bundle 
 $\mu^{-1}_{\frX}(0) \longrightarrow X$ is trivial,
 we have the isomorphism
 \[
  H^n\bigl(\bigoplus_{j \in \Z} H^0(C^{\bullet,j}(\scW_{\frX}(\frU)), d_{-}), d_{+}\bigr)
 \simeq \scW_{X,c}(U) \otimes_{\C} H_{DR}^n(G) 
 \]
 where $H_{DR}^{\bullet}(G)$ is the (algebraic) de Rham cohomology of $G$.
\begin{proof}
 By the definition of the differential $d_{+}$ in \refeq{eq:1} and \refprop{prop:1}, the complex
 $(\bigoplus_j H^0(C^{\bullet, j}, d_-), d_+)$ is a Lie algebra cohomology
 complex associated to the $\frg$-action on $\scL_c(\frU)$.
 By the above isomorphism \refeq{eq:2}, we have
 \[
  H^n(\frg, \scL_c(\frU)) \simeq 
 H^n(\frg, \C[G]) \otimes_{\C} \scL_c(\frU)^G \simeq
 H^n(\frg, \C[G]) \otimes_{\C} \scW_{X,c}(U).
 \]
 Since the cochain complex for the Lie algebra cohomology 
 $H^q(\frg, \C[G])$ is the algebraic de Rham complex for $G$, we obtain
 the isomorphism of the proposition.
\end{proof}
\end{proposition}

\begin{theorem}
\label{thm:local-BRST}
 We have the following isomorphism:
 \[
  \calH^n_{BRST,c}(\frg, \scW_{\frX})
 \simeq \scW_{X,c} \otimes_{\C} H_{DR}^n(G).
 \]
\begin{proof}
 Fix an arbitrary point $x \in X$. Since the algebraic group $G$ acts on
 $\mu^{-1}_{\frX}(0)$ freely, we can take an open subset $U \subset X$ and
 $\frU \subset \frX$ satisfying the assumption of \refprop{prop:2}.
 Consider the spectral sequence $E^{p,q}_r$ defined in \refsec{sec:double-complex}.
 Then, by \refprop{prop:2}, we have
 \[
  E^{p,q}_2 \simeq H^{p}(H^q(C(\scW_{\frX}(\frU)), d_{-}), d_{+}) \simeq
 \scW_{X,c}(U) \otimes_{\C} H_{DR}^{p+q}(G)
 \]
By \reflemma{lemma:3}, the spectral sequence
 $E^{p,q}_r$ converges to the BRST cohomology $H^{p+q}_{BRST,c}(\frg, \scW_{\frX}(\frU))$. 
Consider the differential $d^{(2)} : E^{p,q}_2 \longrightarrow E^{p+2,q-1}_2$
of this spectral sequence. By \refprop{prop:1}, we have
$E^{p,q}_2 = 0$ unless $q = 0$. Therefore we have $d^{(2)} = 0$ and 
the spectral sequence $E^{p,q}_r$ collapses at $r=2$. Thus we have
\[
 H^{p+q}_{BRST, c}(\frg, \scW_{\frX}(\frU)) \simeq \bigoplus_{p \in \Z} E^{p,q}_{\infty}
 \simeq \bigoplus_{p \in \Z}
 E^{p,q}_2 \simeq \scW_{X,c}(U) \otimes_{\C} H_{DR}^{p+q}(G).
\]
This induces an isomorphism of stalks
$H^{p+q}_{BRST, c}(\frg, \scW_{\frX, p^{-1}(x)}) \simeq (\scW_{X, c})_x \otimes_{\C} H^{p+q}_{DR}(G)$
 at any point $x \in X$,
 and hence, we have the isomorphism 
 $\calH^{n}_{BRST, c}(\frg, \scW_{\frX}) \simeq \scW_{X,c} \otimes_{\C}H_{DR}^{n}(G)$ as sheaves on $X$.
\end{proof}
\end{theorem}

\begin{remark}
 In the above proof, we use the fact that the $G$-bundle 
 $\mu^{-1}_{\frX}(0) \longrightarrow X$ is locally trivial in order to take
 a suitable $\frU$. Since our algebraic groups $G$ are
 direct products of general linear groups so that the local triviality
 holds for such algebraic groups $G$ even in the Zariski topology.
 When $G$ is not a direct product of general linear groups, we may need to
 consider in complex analytic topology or \'{e}tale topology.
\end{remark}

One can consider the BRST cohomology of the algebra $\frD(V)$ with
respect to the $G$-action and the moment map $\mu_{\frD}$. We say that
the DQ-algebra $\scW_{X,c}$ is $\scW$-affine (cf. \cite{KR}) when we have the following
equivalence of the abelian categories 
\begin{equation}
\label{eq:3}
 \Mod^{good}_{F}(\scW_{X,c}) \isoto \frD(X_0, c)\mmod, \qquad \scM \mapsto \gHom_{\scW_{X,c}}(\scW_{X,c}, \scM)^{\C^*}
\end{equation}
where $\Mod^{good}_{F}(\scW_{X,c})$ is a category of good $\scW_{X,c}$-modules with equivariant
$\C^*$-action (F-action) (cf. \cite[Section 2]{KR}) 
and $\frD(X_0, c)\mmod$ is a category of finitely generated 
$\frD(X_0, c)$-modules. Its quasi-inverse functor is given by $\scN \mapsto \scW_{X,c} \otimes_{\frD(X_0,c)} \scN$. It is known that when the parameter $c$ is generic, the DQ-algebra $\scW_{X,c}$
is $\scW$-affine (see \cite{KR}, \cite{BK}, \cite{BLPW}).

In \cite[Section 5.6]{GL}, derived analogue of the above equivalence was studied (see also
essentially the same result by a different approach in \cite{MN}). First we define
a $\C[[\hbar]]$-algebra $\frD_{\hbar}(X_0,c) = \Gamma(X, \scW_{X,c}(0))$. 
For further discussion, we also define
$\frD_{\hbar}(X_0,c)$-modules
$\frL^{\wedge}_{\frX,c} = \Gamma(X, p_*(\scL_c(0))) \simeq \Gamma(\frX, \scL_c(0))$,
and $\frL^{\wedge}_{T^*V,c} = \Gamma(T^*V, \scL_{T^*V,c}(0))$, where
$\scL_{T^*V,c}(0)$ is a left $\scW_{T^*V}(0)$-module defined by
\[
 \scL_{T^*V,c}(0) = \scW_{T^*V}(0) / \sum_{i=1}^{\Dim \frg} \scW_{T^*V}(-1) (\mu_{\scW}(A_i) + c(A_i)).
\]
Note that we have a natural embedding 
$\frL^{\wedge}_{T^*V,c} \subset \frL^{\wedge}_{\frX,c}$ because
$\scW_{\frX} = \scW_{T^*V} \vert_{\frX}$.
Since $\Gamma(X, \scW_{X,c}[\hbar^{1/2}])^{\C^*} = \frD(X_0,c)$ (\refprop{prop:3}), and
\[
 \scW_{X,c} \simeq_{\C} \calO_X \otimes_{\C} \C((\hbar)) \supset
 \calO_X \otimes_{\C} \C[[\hbar]] \simeq \scW_{X,c}(0)
\]
as sheaves of $\C$-vector spaces (see \refsec{sec:w-algebras}), then we have an isomorphism
of $\C$-algebras
\begin{equation}
\label{eq:4}
 (\C[\hbar^{1/2}, \hbar^{-1/2}] \otimes_{\C[\hbar]} \frD_{\hbar}(X_0,c))^{\C^*}
 \simeq \frD(X_0,c).
\end{equation}
We also have a natural isomorphism 
\begin{equation}
 \label{eq:8}
 (\C[\hbar^{1/2}, \hbar^{-1/2}] \otimes_{\C[\hbar]} \frL^{\wedge}_{T^*V,c})^{\C^*}
  \simeq \frL_c.
\end{equation}

Now assume that $\frD_{\hbar}(X_0, c)$ 
is of finite global dimension. Then, we have the following equivalence of triangulated 
categories:
\begin{equation}
\label{eq:6} 
 D^b(\scW_{X,c}(0)\mmod) \isoto D^b(\frD_{\hbar}(X_0,c)\mmod), \qquad
 \scM \mapsto \RHom{}_{\scW_{X,c}(0)}(\scW_{X,c}(0), \scM)
\end{equation}
where $D^b(\calC)$ is the bounded derived category of an abelian category $\calC$.
Its quasi-inverse functor is given by the tensor product functor 
\begin{equation}
\label{eq:7} 
 D^b(\frD_{\hbar}(X_0,c)\mmod) \isoto D^b(\scW_{X,c}(0)\mmod), \qquad
 M \mapsto \scW_{X,c}(0) \LTensor_{\frD_{\hbar}(X_0,c)} M.
\end{equation}

\begin{remark}
In \cite[Section 5.6]{GL}, such a derived equivalence was studied for the rational 
Cherednik algebra associated with the wreath product group 
$G(\ell, 1, n) = (\Z/\ell\Z) \wr \frS_n$. 
Our algebra $\frD(X_0,c)$ does {\em not} coincide with the
rational Cherednik algebra, but with its {\em spherical subalgebra}. However, one can 
easily check that their proof works also for the spherical subalgebra by replacing
their quantized Procesi bundle by $\scW_{X,c}(0)$ when the spherical subalgebra has
finite global dimension. 
Moreover, it is also easy to check that their proof works not only in the case of the rational 
Cherednik algebra for $G(\ell,1,n)$, but also for every other algebras obtained by
quantum Hamiltonian reduction under the basic conditions which we assumed in 
\refsec{sec:hamilt-reduct} if $\frD(X_0,c)$ has finite global dimension.
\end{remark}

Let $c_0$ be a parameter such that the DQ-algebra $\scW_{X,c_0}$ is $\scW$-affine. 
By \cite[Theorem 2.9]{KR}, we have
$R^n\Gamma(X, p_{*}(\scL_{c_0})) = 0$ for $n \ne 0$. Since we have an isomorphism
$p_*(\scL_{c_0}) \simeq p_*(\calO_{\mu_{\frX}^{-1}(0)}) \otimes_{\C} \C((\hbar))$ as sheaves on $X$, it 
induces the following vanishing of cohomologies,
\begin{equation}
\label{eq:5}
 R^n\Gamma(X, p_{*}(\calO_{\mu_{\frX}^{-1}(0)})) = 0 \qquad \text{for } n \ne 0.
\end{equation}

Now let $c$ be an arbitrary parameter such that $\frD(X_0,c)$ has finite global dimension.

\begin{lemma}
 \label{lemma:11}
 We have the isomorphism
 \[
  \RHom{}_{\scW_{X,c}(0)}(\scW_{X,c}(0), p_{*}(\scL_c(0))) \simeq \frL_{\frX,c}^{\wedge}
 \]
 in the derived category $D^b(\frD_{\hbar}(X_0,c)\mmod)$.
\begin{proof}
 Since $p_{*}(\scL_c(0))$ is a locally free $\scW_{X,c}(0)$-module, we have
\begin{align*}
  \RHom{}_{\scW_{X,c}(0)}(\scW_{X,c}(0), p_{*}(\scL_c(0))) 
 &\simeq R\Gamma(X, \lHom_{\scW_{X,c}(0)}(\scW_{X,c}(0), p_{*}(\scL_c(0)))), \\
 &\simeq R\Gamma(X, p_{*}(\scL_c(0))).
\end{align*}
By \cite[Lemma 2.12]{KR}, \refeq{eq:5} induces $R^n\Gamma(X, p_{*}(\scL_c(0))) = 0$ for
$n \ne 0$. Thus $\RHom{}_{\scW_{X,c}(0)}(\scW_{X,c}(0), p_*(\scL_c(0)))$ is isomorphic
to its zeroth cohomology 
\[
 \gHom_{\scW_{X,c}(0)}(\scW_{X,c}(0), p_{*}(\scL_c(0))) \simeq \Gamma(X, p_{*}(\scL_c(0))) \simeq \frL^{\wedge}_{\frX,c}.
\]
\end{proof}
\end{lemma}

\begin{theorem}
\label{thm:global-BRST}
 Assume that $\frD(X_0,c)$ has finite global dimension.
 Then we have the isomorphism of $\C$-algebras
\[
 H^{n}_{BRST, c}(\frg, \frD(V)) \simeq \frD(X_0, c) \otimes_{\C} H_{DR}^n(G).
\]
\begin{proof}
 Consider the derived equivalence \refeq{eq:6}. By \reflemma{lemma:11}, we have
 isomorphisms
 \[
  \scW_{X,c}(0) \LTensor_{\frD_{\hbar}(X_0,c)} \frL^{\wedge}_{\frX,c} \simeq 
 \scW_{X,c}(0) \otimes_{\frD_{\hbar}(X_0,c)} \frL^{\wedge}_{\frX,c} \simeq
 p_*(\scL_c(0))
 \]
 in the derived category $D^b(\scW_{X,c}(0)\mmod)$. Using these isomorphisms
for the Lie algebra cohomology $H^n(\frg, p_*(\scL_c(0)))$, we have isomorphisms
\begin{align*}
 H^n(\frg, \scL_c(0)) & \simeq H^n(\frg, \scW_{X,c}(0) \LTensor_{\frD_\hbar(X_0,c)} 
 \frL^{\wedge}_{\frX,c}), \\
 &\simeq H^n(\frg, \scW_{X,c}(0) \otimes_{\frD_\hbar(X_0,c)} \frL^{\wedge}_{\frX,c}), \\
 &\simeq \scW_{X,c}(0) \otimes_{\frD_\hbar(X_0,c)} H^n(\frg, \frL^{\wedge}_{\frX,c}).
\end{align*}
Here the last isomorphism follows from the fact which $\frg$ acts trivially on $\scW_{X,c}(0)$.
On the other hand, by the proof of \refprop{prop:2}, we have an isomorphism
\[
 H^n(\frg, p_*(\scL_c(0))) \simeq \scW_{X,c}(0) \otimes_{\frD_\hbar(X_0,c)} \frD_\hbar(X_0,c)
 \otimes_{\C} H^n_{DR}(G).
\]
Applying the equivalence \refeq{eq:6} again, we have the following isomorphism
\[
 H^n(\frg, \frL^{\wedge}_{\frX,c}) \simeq \frD_\hbar(X_0,c) \otimes_{\C} H^n_{DR}(G)
 \qquad \text{for } n \in \Z.
\]
We compare Lie algebra cohomologies $H^n(\frg, \frL^{\wedge}_{\frX,c})$ and
$H^n(\frg, \frL^{\wedge}_{T^*V,c})$. Since we assume that the moment map $\mu_{T^*V}$ is flat,
by \cite[Proposition 2.6]{BK} (or equivalently by \refprop{prop:1}), we have
isomorphisms of vector spaces
\[
 \frL^{\wedge}_{T^*V,c} \simeq \bigoplus_{m \ge 0} \C[\mu_{T^*V}^{-1}(0)] \hbar^k
 \subset \frL^{\wedge}_{\frX,c} \simeq \bigoplus_{m \ge 0} \C[\mu_{\frX}^{-1}(0)] \hbar^k.
\]
Moreover, since the filtration of $\frL^{\wedge}_{T^*V,c}$ and $\frL^{\wedge}_{\frX,c}$ is
closed under the $G$-actions, the above isomorphisms are isomorphisms of $G$-modules.
Thus, it is enough to compare $H^n(\frg, \C[\mu_{T^*V}^{-1}(0)])$ and
$H^n(\frg, \C[\mu_{\frX}^{-1}(0)])$.

Note that the $G$-actions are induced from the linear action of $G$ on $V$, the $G$-module
$\C[T^*V]$ is a graded $G$-module with respect to the total degree and images of the dual
moment map $\{\mu_{T^*V}^*(A) \;\vert\; A\in \frg\}$ are homogeneous elements. 
Thus, the $G$-module $\C[\mu^{-1}_{T^*V}(0)]$ is again graded as a $G$-module. Moreover,
since the $G$-action and the $\C^*_{diag}$-action on $T^*V$ commute, $\C[\mu^{-1}_{\frX}(0)]$
is spanned by ratios of two homogeneous polynomials: i.e. we have
\[
 \C[\mu_{\frX}^{-1}(0)] =
 \Span_{\C}\left\{ \frac{f}{g} \;\Bigm\vert\; 
 \begin{array}{l}
  f,g \text{ are homogeneous elements in }\C[\mu^{-1}_{T^*V}] \\
 \text{ such  that $g$ has no zero on $\frX$ and $f$, $g$ are coprime}	       
 \end{array}
 \right\}.
\]
For $p$, $q \in \Z_{\ge 0}$, set
\[
 \C[\mu_{\frX}^{-1}(0)]_{p,q} =
 \Span_{\C}\left\{\frac{f}{g} \;\Bigm\vert\; f, g \text{ are as above and } \deg f = p+q, \,
 \deg g = q \right\}.
\]
Then, $\C[\mu_{\frX}^{-1}(0)]_{p,q}$ is closed under the $G$-action and we have a
decomposition as a $G$-module,
$\C[\mu_{\frX}^{-1}(0)] = \bigoplus_{p,q} \C[\mu_{\frX}^{-1}(0)]_{p,q}$. Setting
$\C[\mu_{T^*V}^{-1}(0)]_{p} = \C[\mu_{\frX}^{-1}(0)]_{p,0} \cap \C[\mu_{T^*V}^{-1}(0)]$,
we also have the decomposition $\C[\mu_{T^*V}^{-1}(0)] = \bigoplus_{p} \C[\mu_{T^*V}^{-1}(0)]_{p}$.
Since $\C[\mu_{\frX}^{-1}(0)]_{p,q}$ is a finite dimensional $G$-module and $G$ is a reductive
group, the $G$-modules $\C[\mu_{\frX}^{-1}(0)]$ and $\C[\mu_{T^*V}^{-1}(0)]$ are completely
reducible. 

Since we assume that $X \longrightarrow X_0$ is birational and $X_0$ is normal,
we have $\C[\mu_{\frX}^{-1}(0)]^{G} = \C[\mu_{T^*V}^{-1}(0)]^{G}$ (see \cite[Lemma 3.1]{BK}).
Thus, we have $\C[\mu_{\frX}^{-1}(0)] = \C[\mu_{T^*V}^{-1}(0)] \oplus N$ for a $G$-module
$N$ such that $N$ is a direct sum of irreducible $G$-modules which are non-trivial.
Since $G$ is a reductive group, we have $H^n(\frg, N) = 0$ for any $n$ by 
\cite[Theorem 7.8.9]{W}. Therefore, we conclude that 
$H^n(\frg, \C[\mu_{T^*V}^{-1}(0)]) \simeq H^n(\frg, \C[\mu_{\frX}^{-1}(0)])$ and hence
\[
 H^n(\frg, \frL^{\wedge}_{T^*V,c}) \simeq H^n(\frg, \frL^{\wedge}_{\frX,c})
 \simeq \frD_{\hbar}(X_0, c) \otimes_{\C} H^n_{DR}(G)
\]
for any $n \in \Z$.

Now apply the functor $(\C[\hbar^{1/2}, \hbar^{-1/2}] \otimes_{\C[\hbar]} \blkbar)^{\C^*}$ to the
above isomorphism. Note that $\hbar$ is invariant under the action of $\C^*$ and the 
differential $d_{+}$ of Lie algebra cochain complex \refeq{eq:1} commutes with the 
$\C^*$-action since $\mu_{\scW}(A_i)$ are $\C^*$-invariant elements. Thus, the
functor $(\C[\hbar^{1/2}, \hbar^{-1/2}] \otimes_{\C[\hbar]} \blkbar)^{\C^*}$ commutes with
the Lie algebra cohomology $H^n(\frg, \blkbar)$. Finally, we have an isomorphism \refeq{eq:4}
and \refeq{eq:8}.
Therefore, we conclude that the above isomorphism implies the isomorphism
\begin{equation}
\label{eq:9} 
 H^n(\frg, \frL_{c}) \simeq \frD(X_0,c) \otimes_{\C} H^n_{DR}(G) 
 \qquad \text{for } n \in \Z.
\end{equation}

Consider the double complex $(C^{\bullet,\bullet}(\frD(V)), d_{+}, d_{-})$, which is 
defined in \refsec{sec:double-complex}, for $R = \frD(V)$. 
Since we assume the moment map $\mu_{T^*V}$ is a flat morphism, for $j \in \Z$ we have
 \[
 H^n(C^{\bullet, j}(\frD(V), d_{-})) \simeq H_{-n}(\frg, \frD(V)) \simeq
 \begin{cases}
  \frL_c \otimes \Lambda^j(\frg) & (n = 0), \\
  0 & (n \ne 0),
 \end{cases}
 \]
 by the same argument of \refsec{sec:vanish-neg-cohom}. Combining it with the isomorphism
 \refeq{eq:9}, the second term of the spectral sequence $E_2^{p,q}$ associated with
 the double complex is isomorphic to
 \[
  E^{p,q}_2 \simeq H^p(\frg, H_{-q}(\frg, \frD(V))) \simeq \frD(X_0,c) \otimes_{\C} H^{p}_{DR}(G)
 \]
 if $q=0$, and $0$ otherwise.
 As in the proof of \refthm{thm:local-BRST}, the spectral sequence $E_r^{p,q}$ collapses at
 $r=2$ and it converges to the BRST cohomology $H^{p+q}_{BRST,c}(\frg, \frD(V))$.
\end{proof}
\end{theorem}

\begin{remark}
 This computation of the BRST cohomology is essentially based on 
 techniques used in \cite{GG} and \cite{G} 
 to compute the Lie algebra homology and cohomology 
 for the finite W-algebras.
 On the other hand, in the case of this paper, the group $G$ does not
 act freely on $\mu_{T^*V}^{-1}(0)$ and this makes an obstruction to
 compute the Lie algebra cohomology. To avoid this difficulty, we make
 use of deformation quantization and the two kinds of $\scW$-affinities.
\end{remark}

\subsection{Examples}
\label{sec:examples}

Applying the results of \refthm{thm:local-BRST} and \refthm{thm:global-BRST}, we can determine
explicitly BRST cohomologies for the algebras defined in \refsec{sec:alg-q-Hamilton}.

Consider the exterior algebra $\Lambda(e_1, e_3, \dots, e_{2m-1})$ 
generated by the generators $e_1$, $\dots$, $e_{2m-1}$ where
$e_i$ is regarded as an element of degree $i$. It is a graded algebra.
The de Rham cohomology of the general linear group $GL(\C^m)$ is given as follows.

\begin{lemma}
\label{lemma:10}
We have the following isomorphism of graded algebras
\[
 H_{DR}^{\bullet}(GL(\C^m)) 
 \simeq \Lambda(e_1, e_3, \dots, e_{2m-1}).
\]
\end{lemma}

\subsubsection{BRST cohomology of the deformed preprojective algebras}

For the deformed preprojective algebra associated to the quiver $Q = (I, E)$
defined in \refsec{sec:deform-preproj-alg}, we consider the reductive
algebraic group 
$G = \prod_{i \in I} GL(\C^{\delta_i}) / \C^*_{diag} \simeq \prod_{i \in I \backslash \{0\}} GL(\C^{\delta_i})$.
Here we use the fact $\delta_0 = 1$ for the last isomorphism. 
By using \reflemma{lemma:10} and the K\"{o}nneth formula, 
it is easy to calculate the de Rham cohomology $H_{DR}^{\bullet}(G)$.
Then, we have the following isomorphisms by \refthm{thm:local-BRST} and \refthm{thm:global-BRST}.
 \begin{align*}
  \calH^{\bullet}_{BRST,c}(\frg, \scW_{\frX}) 
 &\simeq \scW_{X,c} \otimes \bigotimes_{i \in I \backslash \{0\}}
  \Lambda(e_1, \dots, e_{2\delta_i - 1}), \\
  H^{\bullet}_{BRST,c}(\frg, \frD(V)) 
 &\simeq \frD(X_0, c) \otimes \bigotimes_{i \in I \backslash \{0\}}
  \Lambda(e_1, \dots, e_{2\delta_i - 1}).
 \end{align*}

In the case that $Q_0$ is the affine Dynkin diagram of type $A_{\ell}^{(1)}$,
$\delta_i = 1$ for all $i \in I \backslash \{0\} = \{1, \dots, \ell\}$, and
\[
  \calH^m_{BRST,c}(\frg, \scW_{\frX}) 
  \simeq \scW_{X,c}^{\oplus \binom{\ell}{m}}
\]
where $\binom{\ell}{m}$ is the binomial coefficient.

In the case that $Q_0$ is the affine Dynkin diagram of type $D_{\ell}^{(1)}$,
we have $G \simeq (\C^*)^3 \times GL(\C^2)^{\times \ell-3}$, and
\[
   \calH^m_{BRST,c}(\frg, \scW_{\frX}) 
 \simeq \scW_{X,c}^{\oplus \nu_n}
\]
where $\nu_m$ is determined by the generating function
$(1+t)^{\ell}(1+t^3)^{\ell-3} = \sum_{m} \nu_m t^m$.

The type $E_{\ell}^{(1)}$ case can be determined similarly.

\subsubsection{BRST cohomology of the symplectic reflection algebras}

For the symplectic reflection algebra associated to the quiver $Q' = (I', E')$
of rank $n$, i.e. with dimension vector $n \delta$, 
we consider the reductive algebraic group $G = \prod_{i \in I'} GL(\C^{n \delta_i})$. Thus, we have the following isomorphisms
\begin{align*}
   \calH^{\bullet}_{BRST,c}(\frg, \scW_{\frX}) 
 &\simeq \scW_{X,c} \otimes \bigotimes_{i \in I'}
 \Lambda(e_1, \dots, e_{2n\delta_i - 1}), \\
   H^{\bullet}_{BRST,c}(\frg, \frD(V)) 
 &\simeq \frD(X_0, c) \otimes \bigotimes_{i \in I'}
 \Lambda(e_1, \dots, e_{2n\delta_i - 1}).
\end{align*} 

For the type $A_{\ell}^{(1)}$ case, $\delta_i = 1$ for all $i$ and we have
\[
   \calH^{m}_{BRST,c}(\frg, \scW_{\frX}) 
 \simeq \scW_{X,c}^{\oplus \nu'_m}
\]
where $\nu'_n$ is determined by the generating function 
$(1+t+t^3+\dots+t^{2n-1})^\ell = \sum_{m} \nu'_m t^m$

\subsubsection{BRST cohomology of the quantization of hypertoric varieties}

Consider the algebra $\frD(X_0, c)$ defined in \refsec{sec:q-hypertoric-var}
associated with a torus action given by the $d \times n$ matrix $M$.
For this algebra, we consider the $d$-dimensional torus $G = (\C^*)^d$.
By \refthm{thm:local-BRST} and \refthm{thm:global-BRST}, we have the isomorphisms
\begin{align*}
 \calH^m_{BRST,c}(\frg, \scW_{\frX}) 
 &\simeq \scW_{X,c}^{\oplus \binom{d}{m}}, \\
 H^m_{BRST,c}(\frg, \frD(V)) 
 &\simeq \frD(X_0,c)^{\oplus \binom{d}{m}}.
\end{align*}

\end{document}